\documentclass[12pt]{amsart}
\input{amssymb.sty}
\title{An example of a solid von Neumann algebra}
\author{Narutaka OZAWA}
\address{Department of Mathematical Sciences,
University of Tokyo, Tokyo 153-8914}
\email{narutaka@ms.u-tokyo.ac.jp}
\date{April 1, 2008}

\newtheorem*{thm*}{Theorem}
\newtheorem*{prop*}{Proposition}
\theoremstyle{definition}
\newtheorem*{defn}{Definition}
\newtheorem*{ack}{Acknowledgment}
\newcommand{\R}{{\mathbb R}}
\newcommand{\RP}{\widehat{\mathbb R}}
\newcommand{\N}{{\mathbb N}}
\newcommand{\Z}{{\mathbb Z}}
\newcommand{\T}{{\mathbb T}}
\newcommand{\PP}{{\mathcal P}}
\newcommand{\LL}{{\mathcal L}}
\newcommand{\G}{\Gamma}

\newcommand{\p}{\varphi}
\newcommand{\SL}{\mathrm{SL}}
\newcommand{\solid}{{\mathcal S}}

\begin{document}
\begin{abstract}
We prove that the group-measure-space von Neumann algebra
$L^\infty(\T^2) \rtimes \SL(2,\Z)$ is solid.
The proof uses topological amenability of the action
of $\SL(2,\Z)$ on the Higson corona of $\Z^2$.
\end{abstract}
\maketitle
\section{Introduction}
Let $\SL(2,\Z)=\{\left[\begin{smallmatrix}a&b\\c&d
\end{smallmatrix}\right]:a,b,c,d\in\Z,\ ad-bc=1\}$
act by linear transformations on the 2-torus $\T^2$ with
the Haar measure, and $L^\infty(\T^2)\rtimes\SL(2,\Z)$ be
the crossed product von Neumann algebra.
Recall that a finite von Neumann algebra is called \emph{solid}
if every diffuse subalgebra has a non-amenable relative commutant.
The main result of this paper is the following, which strengthens
a result in \cite{solid,kurosh}. See \cite{ci} for some application
of this result to ergodic theory.
\begin{thm*}
The von Neumann algebra $L^\infty(\T^2)\rtimes\SL(2,\Z)$ is solid.
\end{thm*}

For the proof of Theorem, we take $L^\infty(\T^2)\rtimes\SL(2,\Z)$
as the group von Neumann algebra of the semidirect product
$\Z^2\rtimes\SL(2,\Z)$ of $\Z^2$ by the linear action of $\SL(2,\Z)$,
and study the behavior of the action at infinity.
This involves the notion of amenability
for a group action on a topological space, which we recall briefly.
We refer the reader to \cite{adr1,adr2,bo} for detailed accounts
of amenable actions.
For a discrete group $\G$, we denote by
\[
\PP(\G)=\{ \mu\in\ell_1(\G) : \mu\geq0,\ \|\mu\|=1\}
\]
the space of probability measures on $\G$, 
equipped with the norm topology (which coincides with
the pointwise-convergence topology).
The group $\G$ acts on $\PP(\G)$ by left translations:
$(g\mu)(h)=\mu(g^{-1}h)$ for $g,h\in\G$ and $\mu\in\PP(\G)$.

\begin{defn}
Let $\G$ be a countable discrete group and $X$ be a compact
topological space on which $\G$ acts as homeomorphisms.
We say the $\G$-action (or the $\G$-space $X$) is
\emph{amenable} if there is a sequence of continuous maps
$\mu_n\colon X\to\PP(\G)$ such that
\[
\forall g\in\G,\
\lim_{n\to\infty}\sup_{x\in X}\|\mu_n(gx)-g\mu_n(x)\|=0.
\]
\end{defn}
We consider the linear action of $\SL(2,\Z)$ on $\Z^2$.
Since the stabilizer subgroups of non-zero elements are all cyclic
(amenable), it is easy to show the action of $\SL(2,\Z)$ on the
Stone-\v Cech remainder $\beta\Z^2\setminus\Z^2$ of $\Z^2$ is amenable.
We will prove a stronger proposition.
The Higson corona $\partial\Z^2$ is defined to be the maximal quotient
of $\beta\Z^2\setminus\Z^2$, on which $\Z^2$ acts trivially:
\[
C(\partial\Z^2)=\{ f\in\ell_\infty(\Z^2)
 : \forall a\in\Z^2,\ \lim_{x\to\infty}|f(x+a)-f(x)|=0\}/c_0(\Z^2).
\]
The $\SL(2,\Z)$-action on $\Z^2$ naturally gives rise to
an $\SL(2,\Z)$-action on $\partial\Z^2$.
\begin{prop*}
The $\SL(2,\Z)$-action on $\partial\Z^2$ is amenable.
\end{prop*}

\begin{ack}
The essential part of this work was done during the author's
stay at Banff International Research Station for the workshop
``Topics in von Neumann Algebras.'' The author would like to
thank the institute and the organizers for their kind hospitality,
and Professor Asger T\"ornquist for raising a question during the
workshop that dispelled the author's misconception. The travel
was supported by Japan Society for the Promotion of Science.
\end{ack}
\section{Proof of Proposition}
We consider the group $\SL(2,\R)=\{\left[\begin{smallmatrix}a&b\\c&d
\end{smallmatrix}\right] : a,b,c,d\in\R,\ ad-bc=1\}$
acting on the real projective line $\RP=\R\cup\{\infty\}$ by
linear fractional transformations:
\[
\begin{bmatrix}a&b\\c&d\end{bmatrix}\colon t\mapsto \frac{at+b}{ct+d}.
\]
The stabilizer of the point $\infty\in\RP$ is the subgroup $P$
of upper triangular matrices.
Since $P$ is a closed amenable subgroup of $\SL(2,\R)$,
the linear fractional action of $\SL(2,\Z)$
on $\RP\cong\SL(2,\R)/P$ is amenable.
For the proof of this fact, see Example 3.9 in \cite{adr1}
or Section 5.4 in \cite{bo}.
Now, we observe that the map $\p\colon\Z^2\setminus\{0\}\to\RP$,
defined by $\p(\begin{smallmatrix}m\\n\end{smallmatrix})=m/n$,
is $\SL(2,\Z)$-equivariant and satisfies
\[
\lim_{x\to\infty}d(\p(x+a),\p(x))=0
\]
for every $a\in\Z^2$, where $d$ is a fixed metric on $\RP$
which agrees with the topology.
By considering $\p^*\colon C(\RP)\to\ell_\infty(\Z^2)$,
one sees that $\p$ gives rise to an $\SL(2,\Z)$-equivariant
continuous map $\tilde{\p}\colon\partial\Z^2\to\RP$.
It is clear from the definition that amenability of $\RP$
implies that of $\partial\Z^2$.
\hspace*{\fill}$\Box$
\section{Proof of Theorem}
The proof of Theorem is almost a verbatim translation of Section 4
of \cite{kurosh}, and we give it rather sketchily.
For another approach, we refer the reader to Chapter 15 of \cite{bo}.

We follow the notations used in Section 4 of \cite{kurosh} and
plug $C^*_\lambda(\Z^2)$ into $A$ and $\SL(2,\Z)$ into $\G$.
We note that $\G$ is virtually-free and hence $\G\in\solid$, i.e.,
the left-and-right translation action of $\G\times\G$ on the
Stone-\v Cech remainder $\beta\G\setminus\G$ of $\G$ is amenable.
It is proved in \cite{kurosh} that $\G\ltimes\Lambda\in\solid$ if
$\G\in\solid$, $\Lambda$ is amenable,
and there is a map $\zeta\colon\Lambda\to\PP(\G)$ such that
\[
\lim_{y\to\infty}\bigl(\|g\zeta(y)-\zeta(gy)\|
 +\|\zeta(xyx')-\zeta(y)\|\bigr)=0
\]
for all $g\in\G$ and $x,x'\in\Lambda$.
Indeed, for Corollary 4.5 in \cite{kurosh},
the only specific property we require of $\Lambda=\Delta_\G$ is 
the existence of $\xi=\zeta^{1/2}$ in the proof of Proposition 4.4
in \cite{kurosh}. From now on, let $\G=\SL(2,\Z)$ and $\Lambda=\Z^2$
and view them as abstract multiplicative groups.
It is left to construct $\zeta\colon\Lambda\to\PP(\G)$ satisfying
the above condition. Although this can be done by modifying
Proposition 4.1 in \cite{kurosh}, we give an alternative proof here.
By (the proof of) Proposition, there is a sequence of maps
$\zeta_n\colon\Lambda\to\PP(\G)$ such that
\[
\limsup_{y\to\infty}\bigl(\|\zeta_n(gy)-g\zeta_n(y)\|
 +\|\zeta_n(xyx')-\zeta_n(y)\|\bigr)<1/n
\]
for all $n\in\N$, $g\in\G$ and $x,x'\in\Lambda$.
(Indeed, let $\zeta_n(x)=\mu_n(\p(x))$ for a
suitable $\mu_n\colon\RP\to\PP(\SL(2,\Z))$
that verifies amenability of $\RP$.)
For $g\in\G$, $x,x'\in\Lambda$, we define finite subsets
$D_n(g;x,x')\subset\Lambda$ by
\[
D_n(g;x,x')=\{ y\in\Lambda : \|\zeta_n(gy)-g\zeta_n(y)\|
 +\|\zeta_n(xyx')-\zeta_n(y)\|\geq1/n\}.
\]
Take an increasing sequence $\{1\}=E_0\subset E_1\subset\cdots\subset\G$
of finite symmetric subsets such that $\bigcup E_n=\G$ and
likewise for $\{1\}=F_0\subset F_1\subset\cdots\subset\Lambda$.
We define finite subsets $\{1\}=\Omega_0\subset \Omega_1\cdots$ of
$\Lambda$ inductively by
\[
\Omega_n=
\bigcup_{g\in E_n,\,x,x'\in F_n,\,y\in\Omega_{n-1}}
 \bigl(D_n(g;x,x')\cup\{ gy, xyx'\}\bigr)
\]
for $n\geq1$. We define $l(y)=\min\{ n : y\in \Omega_n\}$ and define
$\zeta\colon\Lambda\to\PP(\G)$ by
\[
\zeta(y)=\frac{1}{l(y)}\sum_{n=0}^{l(y)-1}\zeta_n(y).
\]
(The value of $\zeta$ at the unit $1$ does not matter.)
Let $g\in\G$ and $x,x'\in\Lambda$ be given arbitrary and
take $k$ such that $g\in E_k$ and $x,x'\in F_k$.
We observe that $|l(gy)-l(y)|\le1$ and $|l(xyx')-l(y)|\le1$
for every $y$ with $l(y)>k$; and that
$\|\zeta_n(gy)-g\zeta_n(y)\|+\|\zeta_n(xyx')-\zeta_n(y)\|<1/n$
for every $n$ with $k\le n<l(y)$.
It follows that
\[
\lim_{l(y)\to\infty}\bigl(\|g\zeta(y)-\zeta(gy)\|
 +\|\zeta(xyx')-\zeta(y)\|\bigr)=0,
\]
which verifies the required condition.
This proves $\Z^2\rtimes\SL(2,\Z)\in\solid$,
and hence the von Neumann algebra
$L^\infty(\T^2)\rtimes\SL(2,\Z)\cong\LL(\Z^2\rtimes\SL(2,\Z))$
is solid by Theorem 6 in \cite{solid}.
\hspace*{\fill}$\Box$

\end{document}